\documentclass{amsproc}
\usepackage{amssymb}
\usepackage{amsmath}
\usepackage{eurosym}
\usepackage{amsfonts}

            \advance \topmargin by -\headheight \advance
            \topmargin by -\headsep
            \evensidemargin
            \oddsidemargin
\newtheorem{theorem}{Theorem}[section]
\newtheorem{lemma}[theorem]{Lemma}

\newtheorem{corollary}[theorem]{Corollary}

\numberwithin{equation}{section}
\input{tcilatex}

\begin{document}
\title[Gaussian Generalized Tetranacci Numbers]{Gaussian Generalized
Tetranacci Numbers}
\thanks{}
\author[Y\"{u}ksel Soykan]{Y\"{u}ksel Soykan}
\maketitle

\begin{center}
\textsl{Department of Mathematics,}

\textsl{Art and Science Faculty,}

\textsl{Zonguldak B\"{u}lent Ecevit University,}

\textsl{67100, Zonguldak, Turkey}

\textsl{e-mail: \ yuksel\_soykan@hotmail.com}
\end{center}

\textbf{Abstract.} In this paper, we define Gaussian generalized Tetranacci
numbers and as special cases, we investigate Gaussian Tetranacci and
Gaussian Tetranacci-Lucas numbers with their properties.



\section{Introduction and Preliminaries}

In this work, we define Gaussian generalized Tetranacci numbers and give
properties of Gaussian Tetranacci and Gaussian Tetranacci-Lucas numbers as
special cases. First, we present some background about generalized
Tetranacci numbers and Gaussian numbers before defining Gaussian generalized
Tetranacci numbers.

There have been so many studies of the sequences of numbers in the
literature which are defined recursively. Two of these type of sequences are
the sequences of Tetranacci and Tetranacci-Lucas which are special case of
generalized Tetranacci numbers. A generalized Tetranacci sequence $%
\{V_{n}\}_{n\geq 0}=\{V_{n}(V_{0},V_{1},V_{2},V_{3})\}_{n\geq 0}$ is defined
by the fourth-order recurrence relations%
\begin{equation}
V_{n}=V_{n-1}+V_{n-2}+V_{n-3}+V_{n-4},  \label{equation:dcfvtsrewqsxazsae}
\end{equation}%
with the initial values $V_{0}=c_{0},V_{1}=c_{1},V_{2}=c_{2},V_{3}=c_{3}$
not all being zero.

This sequence has been studied by many authors and more detail can be found
in the extensive literature dedicated to these sequences, see for example [%
\ref{bib:hathiwala2017}], [\ref{melham1999}], [\ref{natividad2013}], [\ref%
{bib:singh2014}], [\ref{waddill1967}], [\ref{waddill1992}].

The sequence $\{V_{n}\}_{n\geq 0}$ can be extended to negative subscripts by
defining%
\begin{equation*}
V_{-n}=-V_{-(n-1)}-V_{-(n-2)}-V_{-(n-3)}+V_{-(n-4)}
\end{equation*}%
for $n=1,2,3,....$ Therefore, recurrence (\ref{equation:dcfvtsrewqsxazsae})
holds for all integer $n.$

The first few generalized Tetranacci numbers with positive subscript and
negative subscript are given in the following table:

$%
\begin{array}{cccccccc}
n & 0 & 1 & 2 & 3 & 4 & 5 & ... \\ 
V_{n} & c_{0} & c_{1} & c_{2} & c_{3} & c_{0}+c_{1}+c_{2}+c_{3} & 
c_{0}+2c_{1}+2c_{2}+2c_{3} & ... \\ 
V_{-n} & c_{0} & c_{3}-c_{2}-c_{1}-c_{0} & 2c_{2}-c_{3} & 2c_{1}-c_{2} & 
2c_{0}-c_{1} & 2c_{3}-2c_{2}-2c_{1}-3c_{0} & ...%
\end{array}%
$

We consider two special cases of $V_{n}:$ $V_{n}(0,1,1,2)=M_{n}$ is the
sequence of Tetranacci numbers (sequence A000078 in [\ref{bib:sloane}]) and $%
V_{n}(4,1,3,7)=R_{n}$ is the sequence of Tetranacci-Lucas numbers (A073817
in [\ref{bib:sloane}]). In other words, Tetranacci sequence $%
\{M_{n}\}_{n\geq 0}$ and Tetranacci-Lucas sequence $\{R_{n}\}_{n\geq 0}$ are
defined by the fourth-order recurrence relations%
\begin{equation}
M_{n}=M_{n-1}+M_{n-2}+M_{n-3}+M_{n-4},\text{ \ \ \ \ }%
M_{0}=0,M_{1}=1,M_{2}=1,M_{3}=2  \label{equati:fvcvxghsbnz}
\end{equation}%
and 
\begin{equation}
R_{n}=R_{n-1}+R_{n-2}+R_{n-3}+R_{n-4},\text{ \ \ \ \ }%
R_{0}=4,R_{1}=1,R_{2}=3,R_{3}=7  \label{equati:pazertvbcunsmn}
\end{equation}%
respectively.

Next, we present the first few values of the Tetranacci and Tetranacci-Lucas
numbers with positive and negative subscripts:%
\begin{equation*}
\begin{array}{cccccccccccccccc}
n & 0 & 1 & 2 & 3 & 4 & 5 & 6 & 7 & 8 & 9 & 10 & 11 & 12 & 13 & ... \\ 
M_{n} & 0 & 1 & 1 & 2 & 4 & 8 & 15 & 29 & 56 & 108 & 208 & 401 & 773 & 1490
& ... \\ 
M_{-n} & 0 & 0 & 0 & 1 & -1 & 0 & 0 & 2 & -3 & 1 & 0 & 4 & -8 & 5 & ... \\ 
R_{n} & 4 & 1 & 3 & 7 & 15 & 26 & 51 & 99 & 191 & 367 & 708 & 1365 & 2631 & 
5071 & ... \\ 
R_{-n} & 4 & -1 & -1 & -1 & 7 & -6 & -1 & -1 & 15 & -19 & 4 & -1 & 31 & -53
& ...%
\end{array}%
\end{equation*}

It is well known that for all integers $n,$ usual Tetranaci and
Tetranacci-Lucas numbers can be expressed using Binet's formulas%
\begin{equation*}
M_{n}=\frac{\alpha ^{n+2}}{(\alpha -\beta )(\alpha -\gamma )(\alpha -\delta )%
}+\frac{\beta ^{n+2}}{(\beta -\alpha )(\beta -\gamma )(\beta -\delta )}+%
\frac{\gamma ^{n+2}}{(\gamma -\alpha )(\gamma -\beta )(\gamma -\delta )}+%
\frac{\delta ^{n+2}}{(\delta -\alpha )(\delta -\beta )(\delta -\gamma )}
\end{equation*}%
(see for example [\ref{bib:hathiwala2017}] or [\ref{bib:zaveri2015}])

or%
\begin{equation}
M_{n}=\frac{\alpha -1}{5\alpha -8}\alpha ^{n-1}+\frac{\beta -1}{5\beta -8}%
\beta ^{n-1}+\frac{\gamma -1}{5\gamma -8}\gamma ^{n-1}++\frac{\delta -1}{%
5\delta -8}\delta ^{n-1}  \label{equation:srtynsbfdtsra}
\end{equation}%
(see for example [\ref{dresden2014}])

and%
\begin{equation*}
R_{n}=\alpha ^{n}+\beta ^{n}+\gamma ^{n}+\delta ^{n}
\end{equation*}%
respectively, where $\alpha ,\beta ,\gamma $ and $\delta $ are the roots of
the cubic equation $x^{4}-x^{3}-x^{2}-x-1=0.$ Moreover,%
\begin{eqnarray*}
\alpha &=&\frac{1}{4}+\frac{1}{2}\omega +\frac{1}{2}\sqrt{\frac{11}{4}%
-\omega ^{2}+\frac{13}{4}\omega ^{-1}}, \\
\beta &=&\frac{1}{4}+\frac{1}{2}\omega -\frac{1}{2}\sqrt{\frac{11}{4}-\omega
^{2}+\frac{13}{4}\omega ^{-1}}, \\
\gamma &=&\frac{1}{4}-\frac{1}{2}\omega +\frac{1}{2}\sqrt{\frac{11}{4}%
-\omega ^{2}-\frac{13}{4}\omega ^{-1}}, \\
\delta &=&\frac{1}{4}-\frac{1}{2}\omega -\frac{1}{2}\sqrt{\frac{11}{4}%
-\omega ^{2}-\frac{13}{4}\omega ^{-1}},
\end{eqnarray*}%
where%
\begin{equation*}
\omega =\sqrt{\frac{11}{12}+\left( \frac{-65}{54}+\sqrt{\frac{563}{108}}%
\right) ^{1/3}+\left( \frac{-65}{54}-\sqrt{\frac{563}{108}}\right) ^{1/3}}.
\end{equation*}%
Note that we have the following identities: 
\begin{eqnarray*}
\alpha +\beta +\gamma +\delta &=&1, \\
\alpha \beta +\alpha \gamma +\alpha \delta +\beta \gamma +\beta \delta
+\gamma \delta &=&-1, \\
\alpha \beta \gamma +\alpha \beta \delta +\alpha \gamma \allowbreak \delta
+\beta \gamma \delta &=&1, \\
\alpha \beta \gamma \delta &=&-1.
\end{eqnarray*}

We present an identity related with generalized Tetranacci numbers and
Tetranacci numbers.

\begin{theorem}
\label{theorem:sertyhbvcxszasert}For $n\geq 0$ and $m\geq 0$ the following
identity holds: 
\begin{equation}
V_{m+n}=M_{m-2}V_{n+3}+(M_{m-3}+M_{m-4}+M_{m-5})V_{n+2}+(M_{m-3}+M_{m-4})V_{n+1}+M_{m-3}V_{n}
\label{equati:bvcsdtyagvxdscz}
\end{equation}
\end{theorem}

\textit{Proof.} \ We prove the identity by induction on $m.$ If $m=0$ then%
\begin{equation*}
V_{n}=M_{-2}V_{n+3}+(M_{-3}+M_{-4}+M_{-5})V_{n+2}+(M_{-3}+M_{-4})V_{n+1}+M_{-3}V_{n}
\end{equation*}%
which is true because $M_{-2}=0,$ $M_{-3}=1,$ $M_{-4}=-1,$ $M_{-5}=0$.
Assume that the equaliy holds for all $m\leq k.$ For $m=k+1,$ we have%
\begin{eqnarray*}
V_{(k+1)+n} &=&V_{n+k}+V_{n+k-1}+V_{n+k-2}+V_{n+k-3} \\
&=&(M_{k-2}V_{n+3}+(M_{k-3}+M_{k-4}+M_{k-5})V_{n+2}+(M_{k-3}+M_{k-4})V_{n+1}+M_{k-3}V_{n})
\\
&&+(M_{k-3}V_{n+3}+(M_{k-4}+M_{k-5}+M_{k-6})V_{n+2}+(M_{k-4}+M_{k-5})V_{n+1}+M_{k-4}V_{n})
\\
&&+(M_{k-4}V_{n+3}+(M_{k-5}+M_{k-6}+M_{k-7})V_{n+2}+(M_{k-5}+M_{k-6})V_{n+1}+M_{k-5}V_{n})
\\
&&+(M_{k-5}V_{n+3}+(M_{k-6}+M_{k-7}+M_{k-8})V_{n+2}+(M_{k-6}+M_{k-7})V_{n+1}+M_{k-6}V_{n})
\\
&=&(M_{k-2}+M_{k-3}+M_{k-4}+M_{k-5})V_{n+3} \\
&&+((M_{k-3}+M_{k-4}+M_{k-5}+M_{k-6})+(M_{k-4}+M_{k-5}+M_{k-6}+M_{k-7}) \\
&&\text{ \ \ \ \ \ \ \ \ \ \ \ \ \ \ \ \ \ \ \ \ \ \ \ \ \ \ \ \ \ \ }%
+(M_{k-5}+M_{k-6}+M_{k-7}+M_{k-8}))V_{n+2} \\
&&+((M_{k-3}+M_{k-4}+M_{k-5}+M_{k-6})+(M_{k-4}+M_{k-5}+M_{k-6}+M_{k-7}))V_{n+1}
\\
&&+(M_{k-3}+M_{k-4}+M_{k-5}+M_{k-6})V_{n} \\
&=&M_{k-1}V_{n+3}+(M_{k-2}+M_{k-3}+M_{k-4})V_{n+2}+(M_{k-2}+M_{k-3})V_{n+1}+M_{k-2}V_{n}
\\
&=&M_{(k+1)-2}V_{n+3}+(M_{(k+1)-3}+M_{(k+1)-4}+M_{(k+1)-5})V_{n+2} \\
&&+(M_{(k+1)-3}+M_{(k+1)-4})V_{n+1}+M_{(k+1)-3}V_{n}.
\end{eqnarray*}%
By induction on $m,$ this proves (\ref{equation:cvbustydvazxwq}).

The previous Theorem gives the following results as particular examples: For 
$n\geq 0$ and $m\geq 0,$ we have (taking $V_{n}=M_{n}$) 
\begin{equation*}
M_{m+n}=M_{m-2}M_{n+3}+(M_{m-3}+M_{m-4}+M_{m-5})M_{n+2}+(M_{m-3}+M_{m-4})M_{n+1}+M_{m-3}M_{n}
\end{equation*}%
and (taking $V_{n}=R_{n}$) 
\begin{equation*}
R_{m+n}=M_{m-2}R_{n+3}+(M_{m-3}+M_{m-4}+M_{m-5})R_{n+2}+(M_{m-3}+M_{m-4})R_{n+1}+M_{m-3}R_{n}.
\end{equation*}

Next we present the Binet's formula of the generalized Tetranacci sequence.

\begin{lemma}
\label{lemma:mnbvcdfgsd}The Binet's formula of the generalized Tetranacci
sequence $\{V_{n}\}$ is given as%
\begin{equation*}
V_{n}=M_{n-3}V_{0}+(M_{n-3}+M_{n-4})V_{1}+(M_{n-3}+M_{n-4}+M_{n-5})V_{2}+M_{n-2}V_{3}.
\end{equation*}
\end{lemma}

\textit{Proof.} Take $n=0$ and then replace $n$ with $m$ in Theorem \ref%
{theorem:sertyhbvcxszasert}. 
\endproof%

For another proof of the Lemma \ref{lemma:mnbvcdfgsd}, see [\ref%
{bib:singh2014}]. This Lemma is also a special case of a work on the $n$th $%
k $-generalized Fibonacci number (which is also called $k$-step Fibonacci
number) in [\ref{bacani2015}, Theorem 2.2.].

\begin{corollary}
\label{corollary:fgvbxtrtyuhb}The Binet's formula of the generalized
Tetranacci sequence $\{V_{n}\}$ is given as%
\begin{equation*}
V_{n}=A\alpha ^{n-6}+B\beta ^{n-6}+C\gamma ^{n-6}+D\gamma ^{n-6}
\end{equation*}%
where%
\begin{eqnarray*}
A &=&\frac{\alpha -1}{5\alpha -8}(V_{3}\alpha ^{3}+(V_{0}+V_{1}+V_{2})\alpha
^{2}+(V_{1}+V_{2})\alpha +V_{2}) \\
B &=&\frac{\beta -1}{5\beta -8}(V_{3}\beta ^{3}+(V_{0}+V_{1}+V_{2})\beta
^{2}+(V_{1}+V_{2})\beta +V_{2}) \\
C &=&\frac{\gamma -1}{5\gamma -8}(V_{3}\gamma ^{3}+(V_{0}+V_{1}+V_{2})\gamma
^{2}+(V_{1}+V_{2})\gamma +V_{2}) \\
D &=&\frac{\delta -1}{5\delta -8}(V_{3}\gamma ^{3}+(V_{0}+V_{1}+V_{2})\gamma
^{2}+(V_{1}+V_{2})\gamma +V_{2})
\end{eqnarray*}
\end{corollary}

\textit{Proof.} The proof follows from Lemma \ref{lemma:mnbvcdfgsd}\ and (%
\ref{equation:srtynsbfdtsra}).

In fact, Corollary \ref{corollary:fgvbxtrtyuhb} is a special case of a
result in [\ref{bacani2015}, Remark 2.3.].

Next, we give the ordinary generating function $\sum\limits_{n=0}^{\infty
}a_{n}x^{n}$ of the sequence $V_{n}.$

\begin{lemma}
Suppose that $f_{V_{n}}(x)=\sum\limits_{n=0}^{\infty }a_{n}x^{n}$ is the
ordinary generating function of the generalized Tetranacci sequence $%
\{V_{n}\}_{n\geq 0}.$
Then $f_{V_{n}}(x)$ is given by%
\begin{equation}
f_{V_{n}}(x)=\frac{%
V_{0}+(V_{1}-V_{0})x+(V_{2}-V_{1}-V_{0})x^{2}+(V_{3}-V_{2}-V_{1}-V_{0})x^{3}%
}{1-x-x^{2}-x^{3}-x^{4}}.  \label{equation:mnbvcxersoausnb}
\end{equation}
\end{lemma}

Proof. Using (\ref{equation:dcfvtsrewqsxazsae}) and some calculation, we
obtain 
\begin{equation*}
f_{V_{n}}(x)-xf_{V_{n}}(x)-x^{2}f_{V_{n}}(x)-x^{3}f_{V_{n}}(x)=V_{0}+(V_{1}-V_{0})x+(V_{2}-V_{1}-V_{0})x^{2}
\end{equation*}%
which gives (\ref{equation:mnbvcxersoausnb}).%
\endproof%

The previous Lemma gives the following results as particular examples:
generating function of the Tetranacci sequence $M_{n}$ is%
\begin{equation*}
f_{M_{n}}(x)=\sum_{n=0}^{\infty }M_{n}x^{n}=\frac{x}{1-x-x^{2}-x^{3}-x^{4}}
\end{equation*}%
and generating function of the Tetranacci-Lucas sequence $R_{n}$ is 
\begin{equation*}
f_{R_{n}}(x)=\sum_{n=0}^{\infty }R_{n}x^{n}=\frac{4-3x-2x^{2}-x^{3}}{%
1-x-x^{2}-x^{3}-x^{4}}\text{.}
\end{equation*}

In literature, there have been so many studies of the sequences of Gaussian
numbers. A Gaussian integer $z$ is a complex number whose real and imaginary
parts are both integers, i.e., $z=a+ib,$ $a,b\in \mathbb{Z}$. These numbers
is denoted by $\mathbb{Z}[i]$. The norm of a Gaussian integer $a+ib,$ $%
a,b\in \mathbb{Z}$ is its Euclidean norm, that is, $N(a+ib)=\sqrt{a^{2}+b^{2}%
}=\sqrt{(a+ib)(a-ib)}$. For more information about this kind of integers,
see the work of Fraleigh [\ref{fraleighalgebra1976}].

If we use together sequences of integers defined recursively and Gaussian
type integers, we obtain a new sequences of complex numbers such as Gaussian
Fibonacci, Gaussian Lucas, Gaussian Pell, Gaussian Pell-Lucas and Gaussian
Jacobsthal numbers; Gaussian Padovan and Gaussian Pell-Padovan numbers;
Gaussian Tribonacci numbers.

In 1963, Horadam [\ref{bib:horadam1963aa}] introduced the concept of complex
Fibonacci number called as the Gaussian Fibonacci number. Pethe [\ref%
{pethesome1988}] defined the complex Tribonacci numbers at Gaussian
integers, see also [\ref{gurelphd2015}]. There are other several studies
dedicated to these sequences of Gaussian numbers such as the works in [\ref%
{ascigaussianpoly2013}], [\ref{Berzsenyigauss1977}], [\ref%
{catarinogauspell2018}], [\ref{gurelphd2015}], [\ref{halicigausspell2016}], [%
\ref{halicigaussianpell2018}], [\ref{harmancomplex1981}], [\ref%
{bib:horadam1963aa}], [\ref{jordangauussianfib1965}], [\ref%
{pethegeneralized1986}], [\ref{pethegeneragaus1988}], [\ref%
{bib:soykangaussitri2018}], [\ref{tascigaussianter2017}], [\ref%
{tascigaussianpad2018}], [\ref{yagmurgausmodpel2018}], among others.

\section{Gaussian Generalized Tetranacci Numbers}

Gaussian generalized Tetranacci numbers $\{GV_{n}\}_{n\geq
0}=\{GV_{n}(GV_{0},GV_{1},GV_{2},GV_{3})\}_{n\geq 0}$ are defined by%
\begin{equation}
GV_{n}=GV_{n-1}+GV_{n-2}+GV_{n-3}+GV_{n-4},\text{ \ \ \ }
\label{equati:khcvxzosadreqwavbsxc}
\end{equation}%
with the initial conditions%
\begin{equation*}
\text{\ }%
GV_{0}=c_{0}+i(c_{3}-c_{2}-c_{1}-c_{0}),GV_{1}=c_{1}+ic_{0},GV_{2}=c_{2}+ic_{1},GV_{3}=c_{3}+ic_{2},
\end{equation*}%
not all being zero. The sequences $\{GV_{n}\}_{n\geq 0}$ can be extended to
negative subscripts by defining%
\begin{equation*}
GV_{-n}=-GV_{-(n-1)}-GV_{-(n-2)}-GV_{-(n-3)}+GV_{-(n-4)}
\end{equation*}%
for $n=1,2,3,...$. Therefore, recurrence (\ref{equati:khcvxzosadreqwavbsxc})
hold for all integer $n.$ Note that for $n\geq 0$%
\begin{equation}
GV_{n}=V_{n}+iV_{n-1}.  \label{equation:nvcdfcxdsaxzuysa}
\end{equation}%
and 
\begin{equation*}
GV_{-n}=V_{-n}+iV_{-n-1}
\end{equation*}

The first few generalized Gaussian Tetranacci numbers with positive
subscript and negative subscript are given in the following table:%
\begin{equation*}
\begin{array}{ccc}
n & GV_{n} & GV_{-n} \\ 
0 & c_{0}+i(c_{3}-c_{2}-c_{1}-c_{0}) & c_{0}+i(c_{3}-c_{2}-c_{1}-c_{0}) \\ 
1 & c_{1}+ic_{0} & (c_{3}-c_{2}-c_{1}-c_{0})+i(2c_{2}-c_{3}) \\ 
2 & c_{2}+ic_{1} & 2c_{2}-c_{3}+i(2c_{1}-c_{2}) \\ 
3 & c_{3}+ic_{2} & 2c_{1}-c_{2}+i(2c_{0}-c_{1}) \\ 
4 & \allowbreak c_{0}+c_{1}+c_{2}+c_{3}+ic_{3} & 
2c_{0}-c_{1}+i(2c_{3}-3c_{0}-2c_{1}-2c_{2})%
\end{array}%
\end{equation*}
We consider two special cases of $GV_{n}:$ $GV_{n}(0,1,1+i,2+i)=GM_{n}$ is
the sequence of Gaussian Tetranacci numbers and $%
GV_{n}(4-i,1+4i,3+i,7+3i)=GR_{n}$ is the sequence of Gaussian
Tetranacci-Lucas numbers. We formally define them as follows:

Gaussian Tetranacci numbers are defined by%
\begin{equation}
GM_{n}=GM_{n-1}+GM_{n-2}+GM_{n-3}+GM_{n-4},\text{ \ \ \ }
\end{equation}%
with the initial conditions%
\begin{equation*}
\text{\ }GM_{0}=0,GM_{1}=1,GM_{2}=1+i,GM_{3}=2+i
\end{equation*}%
and Gaussian Tetranacci-Lucas numbers are defined by%
\begin{equation}
GR_{n}=GR_{n-1}+GR_{n-2}+GR_{n-3}+GR_{n-4}\text{\ \ }
\end{equation}%
with the initial conditions%
\begin{equation*}
\text{\ }GR_{0}=4-i,GR_{1}=1+4i,GR_{2}=3+i,GR_{3}=7+3i.
\end{equation*}%
Note that for $n\geq 0$%
\begin{equation*}
GM_{n}=M_{n}+iM_{n-1},\text{ }GR_{n}=R_{n}+iR_{n-1}
\end{equation*}%
and 
\begin{equation*}
GM_{-n}=M_{-n}+iM_{-n-1},GR_{-n}=R_{-n}+iR_{-n-1}.
\end{equation*}%
The first few values of Gaussian Tetranacci numbers with positive and
negative subscript are given in the following table.

$%
\begin{array}{ccccccccccc}
n & 0 & 1 & 2 & 3 & 4 & 5 & 6 & 7 & 8 & 9 \\ 
GM_{n} & 0 & 1 & 1+i & 2+i & 4+2i & 8+4i & 15+8i & 29+15i & 56+29i & 108+56i
\\ 
GM_{-n} & 0 & 0 & i & 1-i & -1 & 0 & 2i & 2-3i & -3+i & 1%
\end{array}%
$

The first few values of Gaussian Tetranacci-Lucas numbers with positive and
negative subscript are given in the following table.

$%
\begin{array}{cccccccccc}
n & 0 & 1 & 2 & 3 & 4 & 5 & 6 & 7 & 8 \\ 
GR_{n} & 4-i & 1+4i & 3+i & 7+3i & 15+7i & 26+15i & 51+26i & 99+51i & 191+99i
\\ 
GR_{-n} & 4-i & -1-i & -1-i & -1+7i & 7-6i & -6-i & -1-i & -1+15i & 15-19i%
\end{array}%
$

The following Theorem presents the generating function of Gaussian
generalized Tetranacci numbers.

\begin{theorem}
The generating function of Gaussian generalized Tetranacci numbers is given
as%
\begin{equation}
f_{GV_{n}}(x)=\frac{%
GV_{0}+(GV_{1}-GV_{0})x+(GV_{2}-GV_{1}-GV_{0})x^{2}+(GV_{3}-GV_{2}-GV_{1}-GV_{0})x^{3}%
}{1-x-x^{2}-x^{3}-x^{4}}.  \label{equation:drtgxczvcdfsmv}
\end{equation}
\end{theorem}

\textit{Proof.} \ Let%
\begin{equation*}
f_{GV_{n}}(x)=\sum_{n=0}^{\infty }GV_{n}x^{n}
\end{equation*}%
be generating function of Gaussian generalized Tetranacci numbers. Then
using the definition of generalized Gaussian Tetranacci numbers, and
substracting $xf(x),$ $x^{2}f(x),$ $x^{3}f(x)$ and $x^{4}f(x)$ from $f(x)$
we obtain (note the shift in the index $n$ in the third line)%
\begin{eqnarray*}
&&(1-x-x^{2}-x^{3}-x^{4})f_{GV_{n}}(x) \\
&=&\sum_{n=0}^{\infty }GV_{n}x^{n}-x\sum_{n=0}^{\infty
}GV_{n}x^{n}-x^{2}\sum_{n=0}^{\infty }GV_{n}x^{n}-x^{3}\sum_{n=0}^{\infty
}GV_{n}x^{n}-x^{4}\sum_{n=0}^{\infty }GV_{n}x^{n} \\
&=&\sum_{n=0}^{\infty }GV_{n}x^{n}-\sum_{n=0}^{\infty
}GV_{n}x^{n+1}-\sum_{n=0}^{\infty }GV_{n}x^{n+2}-\sum_{n=0}^{\infty
}GV_{n}x^{n+3}-\sum_{n=0}^{\infty }GV_{n}x^{n+4} \\
&=&\sum_{n=0}^{\infty }GV_{n}x^{n}-\sum_{n=1}^{\infty
}GV_{n-1}x^{n}-\sum_{n=2}^{\infty }GV_{n-2}x^{n}-\sum_{n=3}^{\infty
}GV_{n-3}x^{n}-\sum_{n=4}^{\infty }GV_{n-4}x^{n} \\
&=&(GV_{0}+GV_{1}x+GV_{2}x^{2}+GV_{3}x^{3})-(GV_{0}x+GV_{1}x^{2}+GV_{2}x^{3})-(GV_{0}x^{2}+GV_{1}x^{3})-GV_{0}x^{3}
\\
&&+\sum_{n=4}^{\infty }(GV_{n}-GV_{n-1}-GV_{n-2}-GV_{n-3}-GV_{n-4})x^{n} \\
&=&GV_{0}+(GV_{1}-GV_{0})x+(GV_{2}-GV_{1}-GV_{0})x^{2}+(GV_{3}-GV_{2}-GV_{1}-GV_{0})x^{3}
\end{eqnarray*}
Rearranging above equation, we get%
\begin{equation*}
f_{GV_{n}}(x)=\frac{%
GV_{0}+(GV_{1}-GV_{0})x+(GV_{2}-GV_{1}-GV_{0})x^{2}+(GV_{3}-GV_{2}-GV_{1}-GV_{0})x^{3}%
}{1-x-x^{2}-x^{3}-x^{4}}.
\end{equation*}%
\endproof%

The previous Theorem gives the following results as particular examples: the
generating function of Gaussian Tetranacci numbers is%
\begin{equation}
f_{GM_{n}}(x)=\frac{x+ix^{2}}{1-x-x^{2}-x^{3}-x^{4}}
\label{equati:csvdxzsdqweazxsaer}
\end{equation}%
and the generating function of Gaussian Tetranacci-Lucas numbers is 
\begin{equation}
f_{GR_{n}}(x)=\frac{-\left( 1+i\right) x^{3}-\left( 2+2i\right) x^{2}-\left(
3-5i\right) x+4-i}{1-x-x^{2}-x^{3}-x^{4}}.
\end{equation}%
The result (\ref{equati:csvdxzsdqweazxsaer}) is already known, (see [\ref%
{tascigaussianter2017}]).

We now present the Binet formula for the Gaussian generalized Tetranacci
numbers.

\begin{theorem}
The Binet formula for the Gaussian generalized Tetranacci numbers is%
\begin{equation*}
GV_{n}=\left( A\alpha ^{n-6}+B\beta ^{n-6}+C\gamma ^{n-6}+D\gamma
^{n-6}\right) +i\left( A\alpha ^{n-6}+B\beta ^{n-6}+C\gamma ^{n-6}+D\gamma
^{n-6}\right)
\end{equation*}%
where $A,B,C$ and $D$ are as in Corollary (\ref{corollary:fgvbxtrtyuhb}).
\end{theorem}

\textit{Proof.} The proof follows from Corollary (\ref%
{corollary:fgvbxtrtyuhb}) and $GV_{n}=V_{n}+iV_{n-1}$.%
\endproof%

The previous Theorem gives the following results as particular examples: the
Binet formula for the Gaussian Tetranacci numbers is%
\begin{equation*}
GM_{n}=\left( 
\begin{array}{c}
\frac{\alpha ^{n+2}}{(\alpha -\beta )(\alpha -\gamma )(\alpha -\delta )}+%
\frac{\beta ^{n+2}}{(\beta -\alpha )(\beta -\gamma )(\beta -\delta )} \\ 
+\frac{\gamma ^{n+2}}{(\gamma -\alpha )(\gamma -\beta )(\gamma -\delta )}+%
\frac{\delta ^{n+2}}{(\delta -\alpha )(\delta -\beta )(\delta -\gamma )}%
\end{array}%
\right) +i\left( 
\begin{array}{c}
\frac{\alpha ^{n+2}}{(\alpha -\beta )(\alpha -\gamma )(\alpha -\delta )}+%
\frac{\beta ^{n+2}}{(\beta -\alpha )(\beta -\gamma )(\beta -\delta )} \\ 
+\frac{\gamma ^{n+2}}{(\gamma -\alpha )(\gamma -\beta )(\gamma -\delta )}+%
\frac{\delta ^{n+2}}{(\delta -\alpha )(\delta -\beta )(\delta -\gamma )}%
\end{array}%
\right)
\end{equation*}%
and the Binet formula for the Gaussian Tetranacci-Lucas numbers is%
\begin{equation*}
GR_{n}=\left( \alpha ^{n}+\beta ^{n}+\gamma ^{n}+\delta ^{n}\right) +i\left(
\alpha ^{n-1}+\beta ^{n-1}+\gamma ^{n-1}+\delta ^{n}\right) .
\end{equation*}

The following Theorem present some formulas of Gaussian generalized
Tetranacci numbers.

\begin{theorem}
\label{theorem:tyhbnvcxdsxaz}For $n\geq 1$ we have the following formulas:

\begin{description}
\item[(a)] (Sum of the Gaussian generalized Tetranacci numbers)%
\begin{equation*}
\sum_{k=1}^{n}GV_{k}=\frac{1}{3}%
(GV_{n+2}+2GV_{n}+GV_{n-1}-GV_{0}+GV_{1}-GV_{3})
\end{equation*}

\item[(b)] $\sum_{k=1}^{n}GV_{2k+1}=\frac{1}{3}%
(2GV_{2n+2}+GV_{2n}-GV_{2n-1}-2GV_{0}-GV_{1}-3GV_{2}+GV_{3})$

\item[(c)] $\sum_{k=1}^{n}GV_{2k}=\frac{1}{3}%
(2GV_{2n+1}+GV_{2n-1}-GV_{2n-2}+GV_{0}-GV_{1}+3GV_{2}-2GV_{3}).$
\end{description}
\end{theorem}

\textit{Proof.} \ 

\begin{description}
\item[(a)] Using the recurrence relation%
\begin{equation*}
GV_{n}=GV_{n-1}+GV_{n-2}+GV_{n-3}+GV_{n-4}
\end{equation*}%
i.e.%
\begin{equation*}
GV_{n-4}=GV_{n}-GV_{n-1}-GV_{n-2}-GV_{n-3}
\end{equation*}%
we obtain%
\begin{eqnarray*}
GV_{0} &=&GV_{4}-GV_{3}-GV_{2}-GV_{1} \\
GV_{1} &=&GV_{5}-GV_{4}-GV_{3}-GV_{2} \\
GV_{2} &=&GV_{6}-GV_{5}-GV_{4}-GV_{3} \\
GV_{3} &=&GV_{7}-GV_{6}-GV_{5}-GV_{4} \\
GV_{4} &=&GV_{8}-GV_{7}-GV_{6}-GV_{5} \\
&&\vdots \\
GV_{n-4} &=&GV_{n}-GV_{n-1}-GV_{n-2}-GV_{n-3} \\
GV_{n-3} &=&GV_{n+1}-GV_{n}-GV_{n-1}-GV_{n-2} \\
GV_{n-2} &=&GV_{n+2}-GV_{n+1}-GV_{n}-GV_{n-1} \\
GV_{n-1} &=&GV_{n+3}-GV_{n+2}-GV_{n+1}-GV_{n} \\
GV_{n} &=&GV_{n+4}-GV_{n+3}-GV_{n+2}-GV_{n+1}.
\end{eqnarray*}%
If we add the equations by side by, we get%
\begin{eqnarray*}
\sum_{k=1}^{n}GV_{k} &=&\frac{1}{3}%
(GV_{n+4}-GV_{n+2}-2GV_{n+1}-GV_{0}+GV_{1}-GV_{3}) \\
&=&\frac{1}{3}(GV_{n+2}+2GV_{n}+GV_{n-1}-GV_{0}+GV_{1}-GV_{3}).
\end{eqnarray*}

\item[(b)] When we use (\ref{equati:khcvxzosadreqwavbsxc}), we obtain the
following equalities:%
\begin{eqnarray*}
GV_{k} &=&GV_{k-1}+GV_{k-2}+GV_{k-3}+GV_{k-4} \\
GV_{4} &=&GV_{3}+GV_{2}+GV_{1}+GV_{0} \\
GV_{6} &=&GV_{5}+GV_{4}+GV_{3}+GV_{2} \\
GV_{8} &=&GV_{7}+GV_{6}+GV_{5}+GV_{4} \\
GV_{10} &=&GV_{9}+GV_{8}+GV_{7}+GV_{6} \\
&&\vdots \\
GV_{2n+2} &=&GV_{2n+1}+GV_{2n}+GV_{2n-1}+GV_{2n-2}.
\end{eqnarray*}%
If we rearrange the above equalities, we obtain%
\begin{eqnarray*}
GV_{3} &=&GV_{4}-GV_{2}-GV_{1}-GV_{0} \\
GV_{5} &=&GV_{6}-GV_{4}-GV_{3}-GV_{2} \\
GV_{7} &=&GV_{8}-GV_{6}-GV_{5}-GV_{4} \\
GV_{9} &=&GV_{10}-GV_{8}-GV_{7}-GV_{6} \\
&&\vdots \\
GV_{2n-1} &=&GV_{2n}-GV_{2n-2}-GV_{2n-3}-GV_{2n-4} \\
GV_{2n+1} &=&GV_{2n+2}-GV_{2n}-GV_{2n-1}-GV_{2n-2}.
\end{eqnarray*}%
Now, if we add the above equations by side by, we get%
\begin{eqnarray*}
\sum_{k=1}^{n}GV_{2k+1} &=&GV_{2n+2}-GV_{2}-\sum_{k=1}^{2n-1}GV_{k}-GV_{0} \\
&=&GV_{2n+2}-GV_{2}-\frac{1}{3}%
(GV_{(2n-1)+4}-GV_{(2n-1)+2}-2GV_{(2n-1)+1}-GV_{0}+GV_{1}-GV_{3})-GV_{0} \\
&=&GV_{2n+2}-GV_{2}-\frac{1}{3}%
(GV_{2n+3}-GV_{2n+1}-2GV_{2n}-GV_{0}+GV_{1}-GV_{3})-GV_{0} \\
&=&-\frac{1}{3}%
(-3GV_{2n+2}+GV_{2n+3}-GV_{2n+1}-2GV_{2n}+2GV_{0}+GV_{1}+3GV_{2}-GV_{3}) \\
&=&\frac{1}{3}%
(3GV_{2n+2}-GV_{2n+3}+GV_{2n+1}+2GV_{2n}-2GV_{0}-GV_{1}-3GV_{2}+GV_{3})
\end{eqnarray*}%
and%
\begin{eqnarray*}
3GV_{2n+2}-GV_{2n+3}+GV_{2n+1}+2GV_{2n}
&=&2GV_{2n+2}+(GV_{2n+2}+GV_{2n+1}+GV_{2n}-GV_{2n+3})+GV_{2n} \\
&=&2GV_{2n+2}+GV_{2n}-GV_{2n-1}
\end{eqnarray*}%
So 
\begin{equation*}
\sum_{k=1}^{n}GV_{2k+1}=\frac{1}{3}%
(2GV_{2n+2}+GV_{2n}-GV_{2n-1}-2GV_{0}-GV_{1}-3GV_{2}+GV_{3}).
\end{equation*}

\item[(c)] Since%
\begin{equation*}
\sum_{k=1}^{n}GV_{2k+1}+\sum_{k=1}^{n}GV_{2k}=\sum_{k=1}^{2n+1}GV_{k}-GV_{1}
\end{equation*}%
we have%
\begin{eqnarray*}
\sum_{k=1}^{n}GV_{k} &=&\frac{1}{3}%
(GV_{n+4}-GV_{n+2}-2GV_{n+1}-GV_{0}+GV_{1}-GV_{3}) \\
&=&\frac{1}{3}(GV_{n+2}+2GV_{n}+GV_{n-1}-GV_{0}+GV_{1}-GV_{3}), \\
\sum_{k=1}^{n}GV_{2k+1} &=&\frac{1}{3}%
(2GV_{2n+2}+GV_{2n}-GV_{2n-1}-2GV_{0}-GV_{1}-3GV_{2}+GV_{3})
\end{eqnarray*}%
\begin{eqnarray*}
\sum_{k=1}^{n}GV_{2k}
&=&\sum_{k=1}^{2n+1}GV_{k}-\sum_{k=1}^{n}GV_{2k+1}-GV_{1} \\
&=&\frac{1}{3}%
(GV_{(2n+1)+4}-GV_{(2n+1)+2}-2GV_{(2n+1)+1}-GV_{0}+GV_{1}-GV_{3}) \\
&&-\frac{1}{3}%
(2GV_{2n+2}+GV_{2n}-GV_{2n-1}-2GV_{0}-GV_{1}-3GV_{2}+GV_{3})-GV_{1} \\
&=&\frac{1}{3}(GV_{2n+5}-GV_{2n+3}-2GV_{2n+2}-GV_{0}+GV_{1}-GV_{3})+\frac{1}{%
3}(-2GV_{2n+2}-GV_{2n} \\
&&+GV_{2n-1}+2GV_{0}+GV_{1}+3GV_{2}-GV_{3}-3GV_{1}) \\
&=&\frac{1}{3}%
(GV_{2n+5}-GV_{2n+3}-2GV_{2n+2}-GV_{0}+GV_{1}-GV_{3}-2GV_{2n+2}-GV_{2n}+GV_{2n-1}
\\
&&+2GV_{0}+GV_{1}+3GV_{2}-GV_{3}-3GV_{1}) \\
&=&\frac{1}{3}%
(GV_{2n+5}-GV_{2n+3}-4GV_{2n+2}-GV_{2n}+GV_{2n-1}+GV_{0}-GV_{1}+3GV_{2}-2GV_{3})
\\
&=&\frac{1}{3}(2GV_{2n+1}+GV_{2n-1}-GV_{2n-2}+GV_{0}-GV_{1}+3GV_{2}-2GV_{3})
\end{eqnarray*}

$\allowbreak $This completes the proof. 
\endproof%
\end{description}

As special cases of above Theorem, we have the following two Corollaries.
First one present some formulas of Gaussian Tetranacci numbers.

\begin{corollary}
For $n\geq 1$ we have the following formulas:

\begin{description}
\item[(a)] (Sum of the Gaussian Tetranacci numbers)%
\begin{equation*}
\sum_{k=1}^{n}GM_{k}=\frac{1}{3}(GM_{n+2}+2GM_{n}+GM_{n-1}-(1+i))
\end{equation*}

\item[(b)] $\sum_{k=1}^{n}GM_{2k+1}=\frac{1}{3}%
(2GM_{2n+2}+GM_{2n}-GM_{2n-1}-2-2i)$

\item[(c)] $\sum_{k=1}^{n}GM_{2k}=\frac{1}{3}%
(2GM_{2n+1}+GM_{2n-1}-GM_{2n-2}-2+i).$
\end{description}
\end{corollary}

Second Corollary gives some formulas of Gaussian Tetranacci-Lucas numbers.

\begin{corollary}
For $n\geq 1$ we have the following formulas:

\begin{description}
\item[(a)] (Sum of the Gaussian Tetranacci-Lucas numbers)%
\begin{equation*}
\sum_{k=1}^{n}GR_{k}=\frac{1}{3}(GR_{n+2}+2GR_{n}+GR_{n-1}-10+2i)
\end{equation*}

\item[(b)] $\sum_{k=1}^{n}GR_{2k+1}=\frac{1}{3}%
(2GR_{2n+2}+GR_{2n}-GR_{2n-1}-11-2i)$

\item[(c)] $\sum_{k=1}^{n}GR_{2k}=\frac{1}{3}%
(2GR_{2n+1}+GR_{2n-1}-GR_{2n-2}-2-8i).$
\end{description}
\end{corollary}

In fact, using the method of the proof of Theorem \ref{theorem:tyhbnvcxdsxaz}%
, we can prove the following formulas of generalized Tetranacci numbers.

\begin{theorem}
For $n\geq 1$ we have the following formulas:

\begin{description}
\item[(a)] (Sum of the generalized Tetranacci numbers)%
\begin{equation*}
\sum_{k=1}^{n}V_{k}=\frac{1}{3}(V_{n+2}+2V_{n}+V_{n-1}-V_{0}+V_{1}-V_{3})
\end{equation*}

\item[(b)] $\sum_{k=1}^{n}V_{2k+1}=\frac{1}{3}%
(2V_{2n+2}+V_{2n}-V_{2n-1}-2V_{0}-V_{1}-3V_{2}+V_{3})$

\item[(c)] $\sum_{k=1}^{n}V_{2k}=\frac{1}{3}%
(2V_{2n+1}+V_{2n-1}-V_{2n-2}+V_{0}-V_{1}+3V_{2}-2V_{3}).$
\end{description}
\end{theorem}

As special cases of above Theorem, we have the following two Corollaries.
First one present some formulas of Tetranacci numbers.

\begin{corollary}
For $n\geq 1$ we have the following formulas:

\begin{description}
\item[(a)] (Sum of the Tetranacci numbers)%
\begin{equation*}
\sum_{k=1}^{n}M_{k}=\frac{1}{3}(M_{n+2}+2M_{n}+M_{n-1}-1)
\end{equation*}

\item[(b)] $\sum_{k=1}^{n}M_{2k+1}=\frac{1}{3}(2M_{2n+2}+M_{2n}-M_{2n-1}-2)$

\item[(c)] $\sum_{k=1}^{n}M_{2k}=\frac{1}{3}(2M_{2n+1}+M_{2n-1}-M_{2n-2}-2).$
\end{description}
\end{corollary}

Second Corollary gives some formulas of Tetranacci-Lucas numbers.

\begin{corollary}
For $n\geq 1$ we have the following formulas:

\begin{description}
\item[(a)] (Sum of the Tetranacci-Lucas numbers)%
\begin{equation*}
\sum_{k=1}^{n}R_{k}=\frac{1}{3}(R_{n+2}+2R_{n}+R_{n-1}-10)
\end{equation*}

\item[(b)] $\sum_{k=1}^{n}R_{2k+1}=\frac{1}{3}(2R_{2n+2}+R_{2n}-R_{2n-1}-11)$

\item[(c)] $\sum_{k=1}^{n}R_{2k}=\frac{1}{3}(2R_{2n+1}+R_{2n-1}-R_{2n-2}-2).$
\end{description}
\end{corollary}

Note that if the sum starts with the zero then the constant in the formula
may only change, for example 
\begin{equation*}
\sum_{k=0}^{n}R_{k}=R_{0}+\sum_{k=1}^{n}R_{k}=4+\frac{1}{3}%
(R_{n+2}+2R_{n}+R_{n-1}-10)=\frac{1}{3}(R_{n+2}+2R_{n}+R_{n-1}+2)
\end{equation*}%
but%
\begin{equation*}
\sum_{k=0}^{n}M_{k}=M_{0}+\sum_{k=1}^{n}M_{k}=\sum_{k=1}^{n}M_{k}=\frac{1}{3}%
(M_{n+2}+2M_{n}+M_{n-1}-1).
\end{equation*}

\section{Some Identities Connecting Gaussian Tetranacci and Gaussian
Tetranacci-Lucas Numbers}

In this section, we obtain some identities of Gaussian Tetranacci numbers
and Gaussian Tetranacci-Lucas numbers.

First, we can give a few basic relations between $\{GM_{n}\}$ and $%
\{GR_{n}\} $\ as 
\begin{equation}
GR_{n}=-GM_{n+3}+6GM_{n+1}-GM_{n}  \label{equation:xserthbcsdfcxzsd}
\end{equation}%
\begin{equation}
GR_{n}=-GM_{n+2}+5GM_{n+1}-2GM_{n}-GM_{n-1}
\label{equation:nbvcfgsvcxuhfdgv}
\end{equation}%
and also%
\begin{equation}
GR_{n}=4GM_{n+1}-3GM_{n}-2GM_{n-1}-GM_{n-2}.
\end{equation}%
Note that the last three identities hold for all integers $n.$ For example,
to show (\ref{equation:xserthbcsdfcxzsd}), writing 
\begin{equation*}
GR_{n}=aGM_{n+3}+bGM_{n+2}+cGM_{n+1}+dGM_{n}
\end{equation*}%
and solving the system of equations%
\begin{eqnarray*}
GR_{0} &=&aGM_{3}+bGM_{2}+cGM_{1}+dGM_{0} \\
GR_{1} &=&aGM_{4}+bGM_{3}+cGM_{2}+dGM_{1} \\
GR_{2} &=&aGM_{5}+bGM_{4}+cGM_{3}+dGM_{2} \\
GR_{3} &=&aGM_{6}+bGM_{5}+cGM_{4}+dGM_{3}
\end{eqnarray*}%
we find that $a=-1,b=0,c=6,d=-1.$ Or using the relations $%
GM_{n}=M_{n}+iM_{n-1}$, $GR_{n}=R_{n}+iR_{n-1}$ and identity $%
R_{n}=-M_{n+3}+6M_{n+1}-M_{n}$ we obtain the identity (\ref%
{equation:xserthbcsdfcxzsd}). The others can be found similarly.

We will present some other identities between Gaussian Tetranacci and
Gaussian Tetranacci-Lucas numbers with the help of generating functions.

The following lemma will help us to derive the generating functions of even
and odd-indexed Gaussian Tetranacci and Gaussian Tetranacci-Lucas sequences.

\begin{lemma}[{[\protect\ref{frontczakconvo2018}]}]
\label{lemma:cvbnosbvxzfczdsaw}Suppose that $f(x)=\sum\limits_{n=0}^{\infty
}a_{n}x^{n}$ is the generating function of the sequence $\{a_{n}\}_{n\geq
0}. $ Then the generating functions of the sequences $\{a_{2n}\}_{n\geq 0}$
and $\{a_{2n+1}\}_{n\geq 0}$ are given as%
\begin{equation*}
f_{a_{2n}}(x)=\sum\limits_{n=0}^{\infty }a_{2n}x^{n}=\frac{f(\sqrt{x})+f(-%
\sqrt{x})}{2}
\end{equation*}%
and%
\begin{equation*}
f_{a_{2n+1}}(x)=\sum\limits_{n=0}^{\infty }a_{2n+1}x^{n}=\frac{f(\sqrt{x}%
)-f(-\sqrt{x})}{2\sqrt{x}}
\end{equation*}%
respectively.
\end{lemma}

The next Theorem presents the generating functions of even and odd-indexed
generalized Tetranacci sequences.

\begin{theorem}
The generating functions of the sequences $V_{2n}$ and $V_{2n+1}$ are given
by%
\begin{equation*}
f_{V_{2n}}(x)=\frac{%
V_{0}+(-3V_{0}+V_{2})x+(-2V_{0}+V_{1}-2V_{2}+V_{3})x^{2}+(-2V_{2}+V_{3})x^{3}%
}{x^{4}+x^{3}-3x^{2}-3x+1}
\end{equation*}%
and%
\begin{equation*}
f_{V_{2n+1}}(x)=\frac{%
V_{1}+(-3V_{1}+V_{3})x+(V_{0}-V_{1}+2V_{2}-V_{3})x^{2}+(V_{0}+V_{1}+V_{2}-V_{3})x^{3}%
}{x^{4}+x^{3}-3x^{2}-3x+1}
\end{equation*}%
respectively.
\end{theorem}

\textit{Proof.} Both statements are consequences of Lemma \ref%
{lemma:cvbnosbvxzfczdsaw} applied to (\ref{equation:mnbvcxersoausnb}) and
some lengthy work.%
\endproof%

From the previous Theorem we get the following results as particular
examples: the generating functions of the sequences $M_{2n}$ and $M_{2n+1}$
are given by%
\begin{equation*}
f_{M_{2n}}(x)=\frac{x^{2}+x}{x^{4}+x^{3}-3x^{2}-3x+1},\text{ }%
f_{M_{2n+1}}(x)=\frac{-x^{2}-x+1}{x^{4}+x^{3}-3x^{2}-3x+1}
\end{equation*}%
and the generating functions of the sequences $R_{2n}$ and $R_{2n+1}$ are
given by%
\begin{equation*}
f_{R_{2n}}(x)=\frac{x^{3}-6x^{2}-9x+4}{x^{4}+x^{3}-3x^{2}-3x+1},\text{ \ }%
f_{R_{2n+1}}(x)=\frac{x^{3}+2x^{2}+4x+1}{x^{4}+x^{3}-3x^{2}-3x+1}.
\end{equation*}

The next Theorem presents the generating functions of even and odd-indexed
Gaussian generalized Tetranacci sequences.

\begin{theorem}
\label{theorem:sgundegtynmscx}The generating functions of the sequences $%
GV_{2n}$ and $GV_{2n+1}$ are given by%
\begin{equation*}
f_{GV_{2n}}=\frac{%
GV_{0}+(-3GV_{0}+GV_{2})x+(-2GV_{0}+GV_{1}-2GV_{2}+GV_{3})x^{2}+(-2GV_{2}+GV_{3})x^{3}%
}{x^{4}+x^{3}-3x^{2}-3x+1}
\end{equation*}%
and%
\begin{equation*}
f_{GV_{2n+1}}=\frac{%
GV_{1}+(-3GV_{1}+GV_{3})x+(GV_{0}-GV_{1}+2GV_{2}-GV_{3})x^{2}+(GV_{0}+GV_{1}+GV_{2}-GV_{3})x^{3}%
}{x^{4}+x^{3}-3x^{2}-3x+1}
\end{equation*}%
respectively.
\end{theorem}

\textit{Proof.} Both statements are consequences of Lemma \ref%
{lemma:cvbnosbvxzfczdsaw} applied to (\ref{equation:drtgxczvcdfsmv}) and
some lengthy algebraic calculations.%
\endproof%

The previous theorem gives the following two corollaries as particular
examples. Firstly, the next one presents the generating functions of even
and odd-indexed Gaussian Tetranacci sequences.

\begin{corollary}
The generating functions of the sequences $GM_{2n}$ and $GM_{2n+1}$ are
given by%
\begin{equation}
f_{GM_{2n}}=\frac{\left( 1+i\right) x+\left( 1-i\right) x^{2}-ix^{3}}{%
x^{4}+x^{3}-3x^{2}-3x+1}  \label{equati:khgfrteosaeczxda}
\end{equation}%
and%
\begin{equation}
f_{GM_{2n+1}}=\frac{1-\left( 1-i\right) x-\left( 1-i\right) x^{2}}{%
x^{4}+x^{3}-3x^{2}-3x+1}  \label{equ:gdservxbzgfsta}
\end{equation}%
respectively.
\end{corollary}

The following Corollary gives the generating functions of even and
odd-indexed Gaussian Tetranacci-Lucas sequences.

\begin{corollary}
The generating functions of the sequences $GR_{2n}$ and $GR_{2n+1}$ are
given by%
\begin{equation}
f_{GR_{2n}}(x)=\frac{(4-i)-\left( 9-4i\right) x-\left( 6-7i\right)
x^{2}+\left( 1+i\right) x^{3}}{x^{4}+x^{3}-3x^{2}-3x+1}
\label{equa:nbvcdfgtsyaazxc}
\end{equation}%
and%
\begin{equation}
f_{GR_{2n+1}}(x)=\frac{(1+4i)+\left( 4-9i\right) x+\left( 2-6i\right)
x^{2}+\left( 1+i\right) x^{3}}{x^{4}+x^{3}-3x^{2}-3x+1}
\label{eq:dfcsxczdrewqazvxb}
\end{equation}%
respectively.
\end{corollary}

The next Corollary present identities between Gaussian Tetranacci and
Gaussian Tetranacci-Lucas sequences.

\begin{corollary}
We have the following identities:%
\begin{eqnarray*}
&&(4-i)GM_{2n}-\left( 9-4i\right) GM_{2n-2}-\left( 6-7i\right)
GM_{2n-4}+\left( 1+i\right) GM_{2n-6} \\
&=&\left( 1+i\right) GR_{2n-2}+\left( 1-i\right) GR_{2n-4}-iGR_{2n-6}, \\
&&(1+4i)GM_{2n}+\left( 4-9i\right) GM_{2n-2}+\left( 2-6i\right)
GM_{2n-4}+\left( 1+i\right) GM_{2n-6} \\
&=&\left( 1+i\right) GR_{2n-1}+\left( 1-i\right) GR_{2n-3}-iGR_{2n-5}), \\
&&(4-i)GM_{2n+1}-\left( 9-4i\right) GM_{2n-1}-\left( 6-7i\right)
GM_{2n-3}+\left( 1+i\right) GM_{2n-5} \\
&=&GR_{2n}-\left( 1-i\right) GR_{2n-2}-\left( 1-i\right) )GR_{2n-4} \\
&&(1+4i)GM_{2n+1}+\left( 4-9i\right) GM_{2n-1}+\left( 2-6i\right)
GM_{2n-3}+\left( 1+i\right) GM_{2n-5} \\
&=&GR_{2n+1}-\left( 1-i\right) GR_{2n-1}-\left( 1-i\right) )GR_{2n-3}
\end{eqnarray*}
\end{corollary}

\textit{Proof.} From (\ref{equati:khgfrteosaeczxda})\ and (\ref%
{equa:nbvcdfgtsyaazxc})\ we obtain%
\begin{equation*}
((4-i)-\left( 9-4i\right) x-\left( 6-7i\right) x^{2}+\left( 1+i\right)
x^{3})f_{GM_{2n}}=(\left( 1+i\right) x+\left( 1-i\right)
x^{2}-ix^{3})f_{GR_{2n}}.
\end{equation*}%
The LHS (left hand side)\ is equal to%
\begin{eqnarray*}
LHS &=&((4-i)-\left( 9-4i\right) x-\left( 6-7i\right) x^{2}+\left(
1+i\right) x^{3})\sum\limits_{n=0}^{\infty }GM_{2n}x^{n} \\
&=&\left( 5-i\right) x^{2}+\left( 5+3i\right) x+\sum\limits_{n=3}^{\infty
}((4-i)GM_{2n}-\left( 9-4i\right) GM_{2n-2} \\
&&-\left( 6-7i\right) GM_{2n-4}+\left( 1+i\right) GM_{2n-6})x^{n}
\end{eqnarray*}%
whereas the RHS is%
\begin{eqnarray*}
RHS &=&(\left( 1+i\right) x+\left( 1-i\right)
x^{2}-ix^{3})\sum\limits_{n=0}^{\infty }GR_{2n}x^{n} \\
&=&\left( 5-i\right) x^{2}+\left( 5+3i\right) x+\sum\limits_{n=3}^{\infty
}(\left( 1+i\right) GR_{2n-2}+\left( 1-i\right) GR_{2n-4}-iGR_{2n-6})x^{n}.
\end{eqnarray*}

Compare the coefficients and the proof of the first identity is done. The
other identities can be proved similarly by using (\ref%
{equati:khgfrteosaeczxda})-(\ref{eq:dfcsxczdrewqazvxb}).%
\endproof%

We present an identity related with Gaussian general Tetranacci numbers and
Tetranacci numbers.

\begin{theorem}
For $n\geq 0$ and $m\geq 0$ the following identity holds: 
\begin{equation}
GV_{m+n}=M_{m-2}GV_{n+3}+(M_{m-3}+M_{m-4}+M_{m-5})GV_{n+2}+(M_{m-3}+M_{m-4})GV_{n+1}+M_{m-3}GV_{n}
\label{equation:cvbustydvazxwq}
\end{equation}
\end{theorem}

\textit{Proof.} \ We prove the identity by strong induction on $m.$ If $m=0$
then%
\begin{equation*}
GV_{n}=M_{-2}GV_{n+3}+(M_{-3}+M_{-4}+M_{-5})GV_{n+2}+(M_{-3}+M_{-4})GV_{n+1}+M_{-3}GV_{n}
\end{equation*}%
which is true because $M_{-2}=0,M_{-3}=1,M_{-4}=-1,M_{-5}=0,$. Assume that
the equaliy holds for $m\leq k.$ For $m=k+1,$ we have%
\begin{eqnarray*}
GV_{(k+1)+n} &=&GV_{n+k}+GV_{n+k-1}+GV_{n+k-2}+GV_{n+k-3} \\
&=&(M_{k-2}GV_{n+3}+(M_{k-3}+M_{k-4}+M_{k-5})GV_{n+2}+(M_{k-3}+M_{k-4})GV_{n+1}+M_{k-3}GV_{n})
\\
&&+(M_{k-3}GV_{n+3}+(M_{k-4}+M_{k-5}+M_{k-6})GV_{n+2}+(M_{k-4}+M_{k-5})GV_{n+1}+M_{k-4}GV_{n})
\\
&&+(M_{k-4}GV_{n+3}+(M_{k-5}+M_{k-6}+M_{k-7})GV_{n+2}+(M_{k-5}+M_{k-6})GV_{n+1}+M_{k-5}GV_{n})
\\
&&+(M_{k-5}GV_{n+3}+(M_{k-6}+M_{k-7}+M_{k-8})GV_{n+2}+(M_{k-6}+M_{k-7})GV_{n+1}+M_{k-6}GV_{n})
\\
&=&(M_{k-2}+M_{k-3}+M_{k-4}+M_{k-5})GV_{n+3} \\
&&+((M_{k-3}+M_{k-4}+M_{k-5}+M_{k-6})+(M_{k-4}+M_{k-5}+M_{k-6}+M_{k-7}) \\
&&+(M_{k-5}+M_{k-6}+M_{k-7}+M_{k-8}))GV_{n+2} \\
&&+((M_{k-3}+M_{k-4}+M_{k-5}+M_{k-6})+(M_{k-4}+M_{k-5}+M_{k-6}+M_{k-7}))GV_{n+1}
\\
&&+(M_{k-3}+M_{k-4}+M_{k-5}+M_{k-6})GV_{n} \\
&=&M_{k-1}GV_{n+3}+(M_{k-2}+M_{k-3}+M_{k-4})GV_{n+2}+(M_{k-2}+M_{k-3})GV_{n+1}+M_{k-2}GV_{n}
\\
&=&M_{(k+1)-2}GV_{n+3}+(M_{(k+1)-3}+M_{(k+1)-4}+M_{(k+1)-5})GV_{n+2} \\
&&+(M_{(k+1)-3}+M_{(k+1)-4})GV_{n+1}+M_{(k+1)-3}GV_{n}
\end{eqnarray*}%
By strong induction on $m,$ this proves (\ref{equation:cvbustydvazxwq}).%
\endproof%

The previous Theorem gives the following results as particular examples: For 
$n\geq 0$ and $m\geq 0,$ we have (taking $GV_{n}=GM_{n}$) 
\begin{equation*}
GM_{m+n}=M_{m-2}GM_{n+3}+(M_{m-3}+M_{m-4}+M_{m-5})GM_{n+2}+(M_{m-3}+M_{m-4})GM_{n+1}+M_{m-3}GM_{n}
\end{equation*}%
and (taking $GV_{n}=GR_{n}$) 
\begin{equation*}
GR_{m+n}=M_{m-2}GR_{n+3}+(M_{m-3}+M_{m-4}+M_{m-5})GR_{n+2}+(M_{m-3}+M_{m-4})GR_{n+1}+M_{m-3}GR_{n}.
\end{equation*}

\section{Matrix Formulation of $V_{n}$}

Now, consider the sequence $\{U_{n}\}$ which is defined by the fourth-order
recurrence relation%
\begin{equation*}
U_{n}=U_{n-1}+U_{n-2}+U_{n-3}+U_{n-4},\text{ \ \ \ \ }%
U_{0}=U_{1}=0,U_{2}=U_{3}=1.
\end{equation*}%
Next, we present the first few values of numbers $U_{n}$ with positive and
negative subscripts:%
\begin{equation*}
\begin{array}{ccccccccccccccccc}
n & 0 & 1 & 2 & 3 & 4 & 5 & 6 & 7 & 8 & 9 & 10 & 11 & 12 & 13 & 14 & ... \\ 
U_{n} & 0 & 0 & 1 & 1 & 2 & 4 & 8 & 15 & 29 & 56 & 108 & 208 & 401 & 773 & 
1490 & ... \\ 
U_{-n} & 0 & 0 & 1 & -1 & 0 & 0 & 2 & -3 & 1 & 0 & 4 & -8 & 5 & -1 & 8 & ...%
\end{array}%
\end{equation*}%
Note that some authors call $\{U_{n}\}$ as a Tetranacci sequence instead of $%
\{M_{n}\}$. The numbers $U_{n}$ can be expressed using Binet's formula%
\begin{equation*}
U_{n}=\frac{\alpha ^{n}}{(\alpha -\beta )(\alpha -\gamma )(\alpha -\delta )}+%
\frac{\beta ^{n}}{(\beta -\alpha )(\beta -\gamma )(\beta -\delta )}+\frac{%
\gamma ^{n}}{(\gamma -\alpha )(\gamma -\beta )(\gamma -\delta )}+\frac{%
\delta ^{n}}{(\delta -\alpha )(\delta -\beta )(\delta -\gamma )}.
\end{equation*}%
The matrix method is very useful method in order to obtain some identities
for special sequences. We define the square matrix $A$ of order $4$ as:%
\begin{equation*}
A=\left( 
\begin{array}{cccc}
1 & 1 & 1 & 1 \\ 
1 & 0 & 0 & 0 \\ 
0 & 1 & 0 & 0 \\ 
0 & 0 & 1 & 0%
\end{array}%
\right) 
\end{equation*}%
such that $\det M=-1.$ Induction proof may be used to establish%
\begin{eqnarray}
A^{n} &=&\left( 
\begin{array}{cccc}
M_{n+1} & M_{n}+M_{n-1}+M_{n-2} & M_{n}+M_{n-1} & M_{n} \\ 
M_{n} & M_{n-1}+M_{n-2}+M_{n-3} & M_{n-1}+M_{n-2} & M_{n-1} \\ 
M_{n-1} & M_{n-2}+M_{n-3}+M_{n-4} & M_{n-2}+M_{n-3} & M_{n-2} \\ 
M_{n-2} & M_{n-3}+M_{n-4}+M_{n-5} & M_{n-3}+M_{n-4} & M_{n-3}%
\end{array}%
\right)   \label{equation:gbhasmnuopdcxsz} \\
&=&\left( 
\begin{array}{cccc}
U_{n+2} & U_{n+1}+U_{n}+U_{n-1} & U_{n+1}+U_{n} & U_{n+1} \\ 
U_{n+1} & U_{n}+U_{n-1}+U_{n-2} & U_{n}+U_{n-1} & U_{n} \\ 
U_{n} & U_{n-1}+U_{n-2}+U_{n-3} & U_{n-1}+U_{n-2} & U_{n-1} \\ 
U_{n-1} & U_{n-2}+U_{n-3}+U_{n-4} & U_{n-2}+U_{n-3} & U_{n-2}%
\end{array}%
\right)  \\
&=&\left( 
\begin{array}{cccc}
U_{n+2} & U_{n+2}-U_{n-2} & U_{n+1}+U_{n} & U_{n+1} \\ 
U_{n+1} & U_{n+1}-U_{n-3} & U_{n}+U_{n-1} & U_{n} \\ 
U_{n} & U_{n}-U_{n-4} & U_{n-1}+U_{n-2} & U_{n-1} \\ 
U_{n-1} & U_{n-1}-U_{n-5} & U_{n-2}+U_{n-3} & U_{n-2}%
\end{array}%
\right) .
\end{eqnarray}

Matrix formulation of $M_{n}$ and $R_{n}$ can be given as%
\begin{equation}
\left( 
\begin{array}{c}
M_{n+3} \\ 
M_{n+2} \\ 
M_{n+1} \\ 
M_{n}%
\end{array}%
\right) =\left( 
\begin{array}{cccc}
1 & 1 & 1 & 1 \\ 
1 & 0 & 0 & 0 \\ 
0 & 1 & 0 & 0 \\ 
0 & 0 & 1 & 0%
\end{array}%
\right) ^{n}\left( 
\begin{array}{c}
M_{3} \\ 
M_{2} \\ 
M_{1} \\ 
M_{0}%
\end{array}%
\right)  \label{equat:nmouyvbcfxdsz}
\end{equation}%
and%
\begin{equation}
\left( 
\begin{array}{c}
R_{n+3} \\ 
R_{n+2} \\ 
R_{n+1} \\ 
R_{n}%
\end{array}%
\right) =\left( 
\begin{array}{cccc}
1 & 1 & 1 & 1 \\ 
1 & 0 & 0 & 0 \\ 
0 & 1 & 0 & 0 \\ 
0 & 0 & 1 & 0%
\end{array}%
\right) ^{n}\left( 
\begin{array}{c}
R_{3} \\ 
R_{2} \\ 
R_{1} \\ 
R_{0}%
\end{array}%
\right) .  \label{equatio:mnbvuhnbtdfx}
\end{equation}%
Induction proofs may be used to establish the matrix formulations $M_{n}$
and $R_{n}$. Note that 
\begin{equation*}
GM_{n}=iU_{n}+U_{n+1}
\end{equation*}%
and%
\begin{equation*}
GR_{n}=\left( 3-2i\right) U_{n+2}-\left( 2-6i\right) U_{n+1}-\left(
1+i\right) U_{n}+\left( 1+i\right) U_{n-2}.
\end{equation*}

Consider the matrices $N_{M},E_{M}$ defined by as follows:%
\begin{eqnarray*}
N_{M} &=&\left( 
\begin{array}{cccc}
2+i & 1+i & 1 & 0 \\ 
1+i & 1 & 0 & 0 \\ 
1 & 0 & 0 & i \\ 
0 & 0 & i & 1-i%
\end{array}%
\right) \\
E_{M} &=&\left( 
\begin{array}{cccc}
GM_{n+3} & GM_{n+2} & GM_{n+1} & GM_{n} \\ 
GM_{n+2} & GM_{n+1} & GM_{n} & GM_{n-1} \\ 
GM_{n+1} & GM_{n} & GM_{n-1} & GM_{n-2} \\ 
GM_{n} & GM_{n-1} & GM_{n-2} & GM_{n-3}%
\end{array}%
\right)
\end{eqnarray*}

Next Theorem presents the relations between $A^{n},N_{M}\ $\ and $E_{M}.$

\begin{theorem}
For $n\geq 3,$ we have%
\begin{equation*}
A^{n}N_{M}=E_{M}.
\end{equation*}
\end{theorem}

\textit{Proof.} \ Using the relation%
\begin{equation*}
GM_{n}=iU_{n}+U_{n+1},
\end{equation*}%
and the calculations 
\begin{eqnarray*}
a &=&\left( 2+i\right) U_{n}+\left( 1+i\right) U_{n-1}+\left( 2+i\right)
U_{n+1}+\left( 2+i\right) U_{n+2} \\
&=&2U_{n}+iU_{n}+U_{n-1}+iU_{n-1}+2U_{n+1}+iU_{n+1}+2U_{n+2}+iU_{n+2} \\
&=&i(U_{n+2}+U_{n+1}+U_{n}+U_{n-1})+(2U_{n+2}+2U_{n+1}+2U_{n}+U_{n-1}) \\
&=&iU_{n+3}+(U_{n+2}+U_{n+1}+U_{n}+(U_{n+2}+U_{n+1}+U_{n}+U_{n-1})) \\
&=&iU_{n+3}+(U_{n+2}+U_{n+1}+U_{n}+U_{n+3})=iU_{n+3}+U_{n+4}=GM_{n+3}
\end{eqnarray*}%
and%
\begin{eqnarray*}
U_{n}+U_{n-1}+U_{n+1}+\left( 1+i\right) U_{n+2}
&=&U_{n}+U_{n-1}+U_{n+1}+U_{n+2}+iU_{n+2} \\
&=&iU_{n+2}+U_{n+3}=GM_{n+2},
\end{eqnarray*}%
we get 
\begin{eqnarray*}
A^{n}N_{M} &=&\left( 
\begin{array}{cccc}
U_{n+2} & U_{n+1}+U_{n}+U_{n-1} & U_{n+1}+U_{n} & U_{n+1} \\ 
U_{n+1} & U_{n}+U_{n-1}+U_{n-2} & U_{n}+U_{n-1} & U_{n} \\ 
U_{n} & U_{n-1}+U_{n-2}+U_{n-3} & U_{n-1}+U_{n-2} & U_{n-1} \\ 
U_{n-1} & U_{n-2}+U_{n-3}+U_{n-4} & U_{n-2}+U_{n-3} & U_{n-2}%
\end{array}%
\right) \left( 
\begin{array}{cccc}
2+i & 1+i & 1 & 0 \\ 
1+i & 1 & 0 & 0 \\ 
1 & 0 & 0 & i \\ 
0 & 0 & i & 1-i%
\end{array}%
\right) \\
&=&\left( 
\begin{array}{cccc}
GM_{n+3} & GM_{n+2} & GM_{n+1} & GM_{n} \\ 
GM_{n+2} & GM_{n+1} & GM_{n} & GM_{n-1} \\ 
GM_{n+1} & GM_{n} & GM_{n-1} & GM_{n-2} \\ 
GM_{n} & GM_{n-1} & GM_{n-2} & GM_{n-3}%
\end{array}%
\right) .
\end{eqnarray*}%
\endproof%

Above Theorem can be proved by mathematical induction as well.

Consider the matrices $N_{R},E_{R}$ defined by as follows:%
\begin{eqnarray*}
N_{R} &=&\left( 
\begin{array}{cccc}
7+3i & 3+i & 1+4i & 4-i \\ 
3+i & 1+4i & 4-i & -1-i \\ 
1+4i & 4-i & -1-i & -1-i \\ 
4-i & -1-i & -1-i & -1+7i%
\end{array}%
\right) \\
E_{R} &=&\left( 
\begin{array}{cccc}
GR_{n+3} & GR_{n+2} & GR_{n+1} & GR_{n} \\ 
GR_{n+2} & GR_{n+1} & GR_{n} & GR_{n-1} \\ 
GR_{n+1} & GR_{n} & GR_{n-1} & GR_{n-2} \\ 
GR_{n} & GR_{n-1} & GR_{n-2} & GR_{n-3}%
\end{array}%
\right) .
\end{eqnarray*}%
The following Theorem presents the relations between $A^{n},$ $N_{R}\ $\ and 
$E_{R}.$

\begin{theorem}
We have%
\begin{equation*}
A^{n}N_{R}=E_{R}.
\end{equation*}
\end{theorem}

\textit{Proof. }The proof requires some lengthy calculation, so we omit it.%
\endproof%

The previous Theorem, also, can be proved by mathematical induction.

Similarly, matrix formulation of $V_{n}$ can be given as%
\begin{equation*}
\left( 
\begin{array}{c}
V_{n+3} \\ 
V_{n+2} \\ 
V_{n+1} \\ 
V_{n}%
\end{array}%
\right) =\left( 
\begin{array}{cccc}
1 & 1 & 1 & 1 \\ 
1 & 0 & 0 & 0 \\ 
0 & 1 & 0 & 0 \\ 
0 & 0 & 1 & 0%
\end{array}%
\right) ^{n}\left( 
\begin{array}{c}
V_{3} \\ 
V_{2} \\ 
V_{1} \\ 
V_{0}%
\end{array}%
\right)
\end{equation*}

Consider the matrices $N_{V},E_{V}$ defined by as follows:%
\begin{eqnarray*}
N_{V} &=&\left( 
\begin{array}{cccc}
{\small ic}_{2}{\small +c}_{3} & {\small ic}_{1}{\small +c}_{2} & {\small ic}%
_{0}{\small +c}_{1} & {\small a}_{1} \\ 
{\small ic}_{1}{\small +c}_{2} & {\small ic}_{0}{\small +c}_{1} & {\small %
(1-i)c}_{0}{\small -ic}_{1}{\small -ic}_{2}{\small +ic}_{3} & {\small a}_{2}
\\ 
{\small ic}_{0}{\small +c}_{1} & {\small (1-i)c}_{0}{\small -ic}_{1}{\small %
-ic}_{2}{\small +ic}_{3} & {\small (1-i)c}_{3}{\small -c}_{1}{\small -(1-2i)c%
}_{2}{\small -c}_{0} & {\small a}_{3} \\ 
{\small (1-i)c}_{0}{\small -ic}_{1}{\small -ic}_{2}{\small +ic}_{3} & 
{\small (1-i)c}_{3}{\small -c}_{1}{\small -(1-2i)c}_{2}{\small -c}_{0} & 
{\small 2ic}_{1}{\small +(2-i)c}_{2}{\small -c}_{3} & {\small a}_{4}%
\end{array}%
\right) , \\
E_{V} &=&\left( 
\begin{array}{cccc}
GV_{n+3} & GV_{n+2} & GV_{n+1} & GV_{n} \\ 
GV_{n+2} & GV_{n+1} & GV_{n} & GV_{n-1} \\ 
GV_{n+1} & GV_{n} & GV_{n-1} & GV_{n-2} \\ 
GV_{n} & GV_{n-1} & GV_{n-2} & GV_{n-3}%
\end{array}%
\right) .
\end{eqnarray*}
where%
\begin{eqnarray*}
a_{1} &=&(1-i)c_{0}-ic_{1}-ic_{2}+ic_{3} \\
a_{2} &=&(1-i)c_{3}-c_{1}-(1-2i)c_{2}-c_{0} \\
a_{3} &=&2ic_{1}+(2-i)c_{2}-c_{3} \\
a_{4} &=&2ic_{0}+(2-i)c_{1}-c_{2}.
\end{eqnarray*}

We now present our final Theorem.

\begin{theorem}
We have%
\begin{equation*}
A^{n}N_{V}=E_{V}.
\end{equation*}
\end{theorem}

\textit{Proof. }The proof requires some lengthy work, so we omit it.%
\endproof%

\end{document}